\documentclass[12pt]{amsart}
\usepackage{amssymb}

\calclayout

\theoremstyle{plain}
\newtheorem{main}{Theorem} 
\newtheorem{mcoro}[main]{Corollary}
\newtheorem{theo}{Theorem}[section]
\newtheorem{prop}[theo]{Proposition}
\newtheorem{lemm}[theo]{Lemma}
\newtheorem{coro}[theo]{Corollary}

\theoremstyle{definition}

\newtheorem{exam}[theo]{Example}

\theoremstyle{remark}
\newtheorem*{rema}{Remark}

\newcommand{\field}[1]{\mathbb{#1}}
\newcommand{\C}{\field{C}}
\newcommand{\R}{\field{R}}

\newcommand{\side}[1]{{}^{#1}\hskip-.075em}

\DeclareMathOperator{\GL}{GL} \DeclareMathOperator{\Hom}{Hom}
\DeclareMathOperator{\Irr}{Irr} \DeclareMathOperator{\Vect}{Vect}
\DeclareMathOperator{\ind}{ind} \DeclareMathOperator{\res}{res}
\DeclareMathOperator{\Rep}{Rep}

\begin{document}

\title{Classification of Equivariant Complex Vector Bundles over a Circle}

\date{October 1, 1999}

\author[J.-H.~Cho]{Jin-Hwan Cho}
\address{Department of Mathematics, Osaka City University}
\email{chofchof@free.kaist.ac.kr}

\author[S.~S.~Kim]{Sung Sook Kim}
\address{Department of Applied Mathematics, Paichai University}
\email{sskim@www.paichai.ac.kr}

\author[M.~Masuda]{Mikiya Masuda}
\address{Department of Mathematics, Osaka City University}
\email{masuda@sci.osaka-cu.ac.jp}

\author[D.~Y.~Suh]{Dong Youp Suh}
\address{Department of Mathematics, Korea Advanced Institute of
Science and Technology}
\email{dysuh@action.kaist.ac.kr}

\thanks{J.-H.~Cho is partially supported by the Korea Science and
Engineering Foundation.  S.~S.~Kim is partially supported by the
Basic Science Research Institute Program, Ministry of Education,
1997 Project number BSRI-97-1428.  D.~Y.~Suh is partially
supported by the Korea Science and Engineering Foundation Grant
971-0103-013-2.}



\begin{abstract}
In this paper we characterize the fiber representations of
equivariant complex vector bundles over a circle and classify
these bundles. We also treat the triviality of equivariant complex
vector bundles over a circle by investigating the extensions of
representations. As a corollary of our results, we calculate the
reduced equivariant $K$-group of a circle for any compact Lie
group.
\end{abstract}

\maketitle

\section{Introduction}

The classification of vector bundles over a topological space is
one of the fundamental problems in topology, and many theories
have been developed to solve the problem.  However the problem
becomes more complex and difficult when one considers it in the
equivariant category.  For instance, every complex vector bundle
over a circle is trivial, but equivariant ones are abundant and
not necessarily trivial.  In this paper, we classify equivariant
complex vector bundles over a circle.  The real case is treated in
another paper~\cite{CKMS99b}.

In order to state our main results, let us fix some notation and
terminology.  Let $G$ be a compact Lie group and let $\rho\colon
G\to O(2)$ be an orthogonal representation of $G$.  The unit
circle of the corresponding $G$-module is denoted by $S(\rho)$. We
set $H=\rho^{-1}(1)$, so that $H$ acts trivially on $S(\rho)$ and
the fiber $H$-module of a complex $G$-vector bundle over $S(\rho)$
is determined uniquely up to isomorphism.  On the other hand, for
a character $\chi$ of $H$ and $g\in G$, a new character
$\side{g}\chi$ of $H$ is defined by
$\side{g}\chi(h)=\chi(g^{-1}hg)$ for $h\in H$.  We say that the
character $\chi$ is \emph{$G$-invariant} if $\side{g}\chi=\chi$
for all $g\in G$.  Our first main theorem characterizes the fiber
$H$-module of a complex $G$-vector bundle over $S(\rho)$.

\begin{main} \label{main:condition_of_fiber_$H$-module}
A complex $H$-module is the fiber $H$-module of a complex
$G$-vector bundle over $S(\rho)$ if and only if its character is
$G$-invariant.
\end{main}

We need more notation to state our second main theorem.  Let
$\Irr(H)$ be the set of characters of irreducible $H$-modules.  It
has a $G$-action defined above.  Since a character is a class
function, the isotropy subgroup $G_\chi$ of $G$ at
$\chi\in\Irr(H)$ contains $H$.  We choose and fix a representative
from each $G$-orbit in $\Irr(H)$ and denote the set of those
representatives by $\Irr(H)/G$.  Denote by $\Vect_G(X)$ the set of
isomorphism classes of complex $G$-vector bundles over a connected
$G$-space $X$ and by $\Vect_{G_{\chi}}(X,\chi)$ the subset of
$\Vect_{G_\chi}(X)$ with a multiple of $\chi$ as the character of
fiber $H$-modules.  They are semi-groups under Whitney sum.  The
decomposition of a $G$-vector bundle into the $\chi$-isotypical
components induces an isomorphism
\[
\Vect_G(X)\cong \prod_{\chi\in\Irr(H)/G}\Vect_{G_{\chi}}(X,\chi).
\]
This reduces the study of $\Vect_G(X)$ to that of
$\Vect_{G_\chi}(X,\chi)$ (see Section~2).

\begin{main} \label{main:semi_group_structure}
The semi-group $\Vect_{G_\chi}(S(\rho),\chi)$ is generated by
\begin{enumerate}
\item one element $L_\chi$ if $\rho(G_\chi)\subset SO(2)$,
\item two elements $L_\chi^\pm$ if $\rho(G_\chi)=O(2)$, and
\item four elements $L_\chi^{\pm\pm}$ with a relation $L_\chi^{++}+
L_\chi^{--}=L_\chi^{+-}+L_\chi^{-+}$ otherwise.
\end{enumerate}
\end{main}

Using this theorem, one can easily enumerate complex $G$-vector
bundles over $S(\rho)$ with a fixed $H$-module as the fiber
$H$-module (see Corollary~\ref{coro:enumeration}).

Our last main theorem is about the triviality of $G$-vector
bundles.  Here a $G$-vector bundle over $X$ is said to be
\emph{trivial} if it is isomorphic to a product bundle with a
$G$-module as its fiber.

\begin{main} \label{main:triviality_for_L}
The triviality of the generators appeared in
Theorem~\ref{main:semi_group_structure} is as follows.
\begin{enumerate}
\item $L_\chi$ is trivial.
\item $L_\chi^{\pm}$ are both trivial or both nontrivial.
\item Two of $L_\chi^{\pm\pm}$ are trivial and the other two are
nontrivial if $|\rho(G_\chi)|/2$ is odd, and $L_\chi^{\pm\pm}$ are
all trivial or all nontrivial if $|\rho(G_\chi)|/2$ is even, where
$|\rho(G_\chi)|$ denotes the order of the dihedral group
$\rho(G_\chi)$.
\end{enumerate}
\end{main}

Since $\rho(G_\chi)\subset SO(2)$ for any $\chi$ if $\rho(G)
\subset SO(2)$, it follows that

\begin{mcoro} \label{main:triviality_SO(2)}
Every $G$-vector bundle over $S(\rho)$ is trivial if
$\rho(G)\subset SO(2)$.
\end{mcoro}

\noindent The reader will find that there are many nontrivial
$G$-vector bundles as well as trivial ones unless $\rho(G)\subset
SO(2)$.

This paper is organized as follows.  Sections~2 and~3 deal with
results on $G$-vector bundles which hold for an arbitrary base
space.  We also recall some results from representation theory,
which turn out to be closely related to the semi-group structure
on $\Vect_G(S(\rho))$ and the triviality of $G$-vector bundles
over $S(\rho)$.  Theorems~\ref{main:condition_of_fiber_$H$-module}
and \ref{main:semi_group_structure} are proved in Sections~4
and~5.  In Section~6 we present another approach to study
$\Vect_G(S(\rho))$.  The triviality of $G$-vector bundles over
$S(\rho)$ is discussed and the proof of
Theorem~\ref{main:triviality_for_L} is given in Section~7.  We
describe $G$-line bundles explicitly in Section~8.  In Section~9
we apply the general results obtained in the previous sections to
the case when $G$ is abelian.  In Section~10 we determine the
reduced equivariant $K$-group of $S(\rho)$, which extends a result
of Y.~Yang \cite{Yan95} for a finite cyclic group to any compact
Lie group $G$.

The subject treated in this paper is classical and the reader may
wonder why we were led to study this subject.  In fact, we were
concerned with what is called the manifold realization problem. It
asks whether a closed smooth $G$-manifold is equivariantly
diffeomorphic to a non-singular real affine $G$-variety.  This
problem was originally considered in the non-equivariant category
by J.~Nash~\cite{Nas52}, and affirmatively solved by
A.~Tognoli~\cite{Tog73}.  Then R.~Palais~\cite{Pal81} considered
the equivariant case above, and some partial affirmative solutions
are obtained, see~\cite{DoMa95} for instance.  It is even
considered to realize smooth $G$-vector bundles over closed smooth
$G$-manifolds by algebraic ones, which is called the bundle
realization problem, and some partial affirmative solutions are
obtained as well, see~\cite{DMS94}.  Apparently the latter problem
is more general than the former, but they are linked.  Namely we
encounter the bundle realization problem to solve the manifold
realization problem.  For instance, we were faced with realizing
real or complex $G$-line bundles over a circle by algebraic ones
to solve the manifold realization problem for two- or
three-dimensional manifolds~\cite{KiMa94,ChSu97}.  This motivated
us to investigate $G$-vector bundles over a circle.  At the
beginning of this research we suspected that the problem might
already be solved, but there is no literature as far as we know.
We hope that it is worth while publishing the results obtained in
this paper in print.

\section{Decomposition of $G$-vector bundles}

Hereafter we omit the adjective ``complex'' for complex vector
bundles and complex modules since we work in the complex category.
Let $G$ be a compact Lie group and let $H$ be a closed normal
subgroup of $G$.  Given a character $\chi$ of $H$ and $g\in G$, a
new character $\side{g}\chi$ of $H$ is defined by
$\side{g}\chi(h)=\chi(g^{-1}hg)$ for $h\in H$.  This defines an
action of $G$ on the set $\Irr(H)$ of characters of irreducible
$H$-modules.  Since a character is a class function, $H$ acts on
$\Irr(H)$ trivially.  Therefore, the isotropy subgroup of $G$ at
$\chi\in \Irr(H)$, denoted by $G_\chi$, contains $H$.  We choose a
representative from each $G$-orbit in $\Irr(H)$ and denote by
$\Irr(H)/G$ the set of those representatives.

Let $X$ be a connected $G$-space on which $H$ acts trivially. Then
all the fibers of a $G$-vector bundle $E$ over $X$ are isomorphic
as $H$-modules.  We call the unique (up to isomorphism) $H$-module
the \emph{fiber $H$-module} of $E$.  As is well-known, $E$
decomposes according to irreducible $H$-modules.  For $\chi\in
\Irr(H)$, we denote by $E(\chi)$ the $\chi$-isotypical component
of $E$, that is, the largest $H$-subbundle of $E$ with a multiple
of $\chi$ as the character of the fiber $H$-module.  Note that
$gE(\chi)$, that is $E(\chi)$ mapped by $g\in G$, is
$\side{g}\chi$-isotypical component of $E$.  This means that
$E(\chi)$ is actually a $G_\chi$-vector bundle and that
$\bigoplus_{\lambda\in G(\chi)}E(\lambda)$, where $G(\chi)$
denotes the $G$-orbit of $\chi$, is a $G$-subbundle of $E$.  Since
$\bigoplus_{\lambda\in G(\chi)}E(\lambda)$ is nothing but the
induced $G$-vector bundle $\ind_{G_\chi}^GE(\chi)$, we have the
following decomposition
\[
E=\bigoplus_{\chi\in\Irr(H)/G}\ind_{G_{\chi}}^G E(\chi)
\]
as $G$-vector bundles.

\begin{lemm} \label{lemm:induced_decomposition}
Two $G$-vector bundles $E$ and $E'$ over $X$ are isomorphic if and
only if $E(\chi)$ and $E'(\chi)$ are isomorphic as
$G_{\chi}$-vector bundles for each $\chi\in\Irr(H)/G$.  In
particular, $E$ is trivial if and only if $E(\chi)$ is trivial for
each $\chi\in\Irr(H)/G$.
\end{lemm}

\begin{proof}
The necessity is obvious since a $G$-vector bundle isomorphism
$E\to E'$ restricts to a $G_{\chi}$-vector bundle isomorphism
$E(\chi)\to E'(\chi)$, and the sufficiency follows from the fact
that $\ind$ is functorial.
\end{proof}

The observation above can be restated as follows.  Denote by
$\Vect_G(X)$ the set of isomorphism classes of $G$-vector bundles
over $X$, and by $\Vect_{G_{\chi}}(X,\chi)$ the subset of
$\Vect_{G_\chi}(X)$ with a multiple of $\chi$ as the character of
fiber $H$-modules.  They are semi-groups under Whitney sum.  Then
the map sending $E$ to $\prod_{\chi\in\Irr(H)/G} E(\chi)$ gives a
semi-group isomorphism
\[
\Phi\colon \Vect_G(X)\to \prod_{\chi\in\Irr(H)/G}
\Vect_{G_{\chi}}(X,\chi).
\]
This reduces the study of $\Vect_G(X)$ to that of
$\Vect_{G_\chi}(X,\chi)$.

\begin{lemm} \label{lemm:hom_argument}
If there is a $G_{\chi}$-vector bundle over $X$ with $\chi$ as the
character of the fiber $H$-module, then $\Vect_{G_{\chi}}(X,\chi)$
is isomorphic to $\Vect_{G_{\chi}/H}(X)$ as semi-groups.  In fact,
if $L$ is such a $G_{\chi}$-vector bundle, then the map
\[
\Vect_{G_{\chi}}(X,\chi)\to\Vect_{G_{\chi}/H}(X)
\]
sending $E$ to $\Hom_H(L,E)$ gives an isomorphism.
\end{lemm}

\begin{proof}
It is easy to check that the map
\[
\Vect_{G_{\chi}/H}(X)\to\Vect_{G_{\chi}}(X,\chi)
\]
sending $F\in \Vect_{G_{\chi}/H}(X)$ to $L\otimes F$ gives the
inverse of the map in the lemma.
\end{proof}

\begin{rema}
The lemma above does not hold in the real category in general, but
it does if $\chi$ is of real type, i.e., if $\chi$ is the
character of a real irreducible $H$-module with the endomorphism
algebra isomorphic to $\R$.
\end{rema}

We conclude this section with the following well-known fact.

\begin{prop} \label{prop:infinite_case}
If the $G$-action on $X$ is transitive, i.e., $X$ is homeomorphic
to $G/K$ for a closed subgroup $K$, then any $G$-vector bundle
over $X$ is of the form
\[
G\times_K W \to G/K=X
\]
for some $K$-module $W$.  In fact, $W$ is the fiber over a point
of $X$ with $K$ as the isotropy subgroup.
\end{prop}

The proposition above implies that there is an isomorphism
\[
\Vect_G(G/K)\cong \Vect_K(*),
\]
where $*$ denotes the one-point $K$-space.

\section{Fiber $H$-modules and extension}

We say that a character $\chi$ of $H$ is $G$-invariant if
$\side{g}\chi=\chi$ for all $g\in G$.  The following proposition
gives a necessary condition for an $H$-module to be the fiber
$H$-module of a $G$-vector bundle $E$ over a connected $G$-space
$X$.

\begin{prop} \label{prop:$G$-invariance_of_fiber_$H$-module}
The character of the fiber $H$-module of $E$ is $G$-invariant.
\end{prop}

\begin{proof}
Let $x\in X$ and $g\in G$.  We know that the fibers $E_x$ and
$E_{gx}$ of $E$ at $x$ and $gx$ are isomorphic as $H$-modules.  On
the other hand, since $g(g^{-1}hg)v=hgv$ for $v\in E_x$, the map
$g\colon E_x\to E_{gx}$ becomes an $H$-equivariant isomorphism if
we consider an $H$-action on $E_{x}$ given through an automorphism
of $H$ given by $h\to g^{-1}hg$.  This implies the proposition.
\end{proof}

For a group $K$ containing $H$, we say that an $H$-module $V$
\emph{extends} to a $K$-module $W$ (or $W$ is a
\emph{$K$-extension} of $V$) if the restriction $\res_H W$ of $W$
to $H$ is isomorphic to $V$.

There are two reasons why we are concerned with the extension of
an $H$-module.  One is that if $V$ is the fiber $H$-module of $E$
and there is a point in the base space with the isotropy subgroup
$K$ larger than $H$, then the fiber over the point gives a
$K$-extension of $V$.  The other is that if $E$ is trivial, i.e,
isomorphic to a product bundle with a $G$-module as its fiber,
then the fiber $H$-module of $E$ must extend to the $G$-module.

If $H$ is a normal subgroup of $K$ and an $H$-module has a
$K$-extension, then its character must be $K$-invariant because a
character is a class function.  But the converse does not hold in
general.  The following proposition gives an answer to the
converse problem when $K/H$ is isomorphic to a subgroup of $O(2)$.

\begin{prop} \label{prop:extension}
Let $H$ be a normal subgroup of $K$ and let $V$ be an irreducible
$H$-module with $K$-invariant character.
\begin{enumerate}
\item If $K/H$ is finite cyclic or isomorphic to $SO(2)$, then
$V$ has a $K$-extension.
\item If $K/H$ is isomorphic to $O(2)$, then $V$ has two
$K$-extensions or none.
\item If $K/H$ is a dihedral group $D_n$ of order $2n$, then
\begin{enumerate}
\item in case $n$ is odd, $V$ has two $K$-extensions,
\item in case $n$ is even, $V$ has either four $K$-extensions
or none.
\end{enumerate}
\end{enumerate}
In any case, if $W$ is a $K$-extension of $V$, then any
$K$-extension of $V$ is of the form $W\otimes U$, where $U$ is a
one-dimensional $K/H$-module viewed as a $K$-module through the
projection $K\to K/H$, and different $U$'s produce different
$K$-extensions.  Therefore, if $V$ has a $K$-extension, then the
number of $K$-extensions of $V$ agrees with the number of
one-dimensional $K/H$-modules.
\end{prop}

\begin{rema}
(1) The character of an $H$-module is $K$-invariant whenever $K/H$
is isomorphic to $SO(2)$ because $SO(2)$ is connected.  So the
$K$-invariance of the character of $V$ in the proposition above is
unnecessary in this case, i.e., any $H$-module extends to a
$K$-module whenever $K/H$ is isomorphic to $SO(2)$.

(2) The proposition above does not hold in the real category but
it does when $K/H$ is cyclic and of odd order.
\end{rema}

\begin{proof}
We shall prove the latter statement in the proposition.  In this
proof, we do not need the assumption that $K/H$ is isomorphic to a
subgroup of $O(2)$.  Suppose that $V$ has $K$-extensions $W$ and
$W'$.  Then $\Hom_H(W,W')$ is a $K/H$-module and one-dimensional
by Schur's lemma because $\res_HW=\res_HW'=V$ is irreducible. One
can easily check that the map
\[
W\otimes \Hom_H(W,W')\to W'
\]
sending $w\otimes f$ to $f(w)$ is a $K$-linear isomorphism.
Therefore, any $K$-extension of $V$ is the tensor product of $W$
and a one-dimensional $K/H$-module viewed as a $K$-module.
Moreover, since $\Hom_H(W,W\otimes U) \cong U$ as $K/H$-modules,
different one-dimensional $K/H$-modules produce different
$K$-extensions of $V$.  This proves the latter statement.

It is well-known and easy to see that there are two $O(2)$- or
$D_n$-modules of dimension one for $n$ odd, and four $D_n$-modules
of dimension one for $n$ even.  This proves the statements on the
number of $K$-extensions in (2) and (3).

It remains to see that $V$ has a $K$-extension in the cases (1)
and~(3--i).  The case (1) is well-known if $K/H$ is finite cyclic,
see for instance~\cite{Isa76}.  The book~\cite{Isa76} treats only
finite groups (so that $H$ is finite), but some direct proofs work
even if $H$ is infinite.  The reader can find one of the direct
proofs in~\cite{CKS99}.  We suspect that the case (1) is known
even when $K/H$ is isomorphic to $SO(2)$, but there is no
literature as far as we know.  Since the proof we found is rather
long and has independent interest, we will give it
in~\cite{CKS99}.  The case~(3--i) is also well-known,
see~\cite{Isa76} or~\cite{CKS99} for an elementary proof.
\end{proof}

It turns out that the above facts on representation theory greatly
influence the semi-group structure on $\Vect_G(X)$ when $X$ is a
circle with $G$-action.

\section{Fiber $H$-modules of $G$-vector bundles over a circle}

Henceforth, we restrict our concern to $G$-vector bundles over a
circle with $G$-action.  For an orthogonal representation
$\rho\colon G\to O(2)$ of $G$, we denote by $S(\rho)$ the unit
circle of the corresponding representation space.  It is
well-known that a circle with continuous (resp.  smooth)
$G$-action is equivariantly homeomorphic (resp.  diffeomorphic) to
$S(\rho)$ for some representation $\rho$~\cite{Sch84}.  We
identify $S(\rho)$ with the unit circle of the complex plane $\C$,
and denote a point in $S(\rho)$ by $z\in\C$ with absolute value
$1$. Set $H=\rho^{-1}(1)$.

Let us observe the subgroup $\rho(G)$ of $O(2)$.  If $\rho(G)$ is
infinite, then it is either $SO(2)$ or $O(2)$ itself.  Otherwise
$\rho(G)$ is a finite cyclic or dihedral group.  Suppose $\rho(G)$
is a finite cyclic subgroup $Z_n$ of $SO(2)$ of order $n\geq1$.
Choose and fix an element $a\in G$ such that $\rho(a)$ is the
rotation through an angle $2\pi/n$.  Then $G$ is generated by $H$
and $a$ under the relation $a^n\in H$, and all the isotropy
subgroups $G_z$ at $z\in S(\rho)$ are equal to $H$.  If $\rho(G)$
is a proper subgroup of $O(2)$ not contained in $SO(2)$, then we
may assume that $\rho(G)$ is a dihedral subgroup $D_n$ of $O(2)$
generated by the reflection matrix about the $x$-axis and the
rotation matrix through an angle $2\pi/n$.  Choose and fix one
more element $b\in G$ such that $\rho(b)$ is the reflection matrix
about the $x$-axis.  Then $G$ is generated by $H$, $a$ and $b$
under the relations, $a^n$, $b^2$, and $(ab)^2\in H$.  The
isotropy subgroup $G_1$ at $1\in S(\rho)$ is generated by $H$ and
$b$, and $G_\mu$ at $\mu=e^{\pi i/n}\in S(\rho)$ is generated by
$H$ and $ab$.

Here is a characterization of the fiber $H$-modules of $G$-vector
bundles over $S(\rho)$.

\setcounter{main}{0}
\begin{main}
An $H$-module is the fiber $H$-module of a $G$-vector bundle over
$S(\rho)$ if and only if its character is $G$-invariant.
\end{main}

\begin{proof}
The necessity follows from
Proposition~\ref{prop:$G$-invariance_of_fiber_$H$-module}, so we
prove the sufficiency.  Let $V$ be an $H$-module with
$G$-invariant character.  We distinguish three cases according to
$\rho(G)$.

\emph{Case 1:} The case where $\rho(G)$ is infinite, i.e.,
$\rho(G)= SO(2)$ or $O(2)$.  In this case the $G$-action on
$S(\rho)$ is transitive, and the isotropy subgroup $K$ of a point
in $S(\rho)$ is $H$ when $\rho(G)=SO(2)$, and contains $H$ as an
index two subgroup when $\rho(G)=O(2)$.  Therefore $V$ has a
$K$-extension by Proposition~\ref{prop:extension} since the
character of $V$ is $G$-invariant (in particular, $K$-invariant)
and $K/H$ is trivial or of order two.  This together with
Proposition~\ref{prop:infinite_case} implies the existence of a
$G$-vector bundle over $S(\rho)$ with $V$ as the fiber $H$-module.

\emph{Case 2:} The case where $\rho(G)=Z_n$.  In this case
$H=G_1=G_\mu$.  Set $W=\ind_H^G V$ and consider the product
$G$-vector bundle $\underline{W}=S(\rho)\times W$.  Choose an
$H$-submodule isomorphic to $V$ in the fiber of $\underline{W}$ at
$1$, and identify it with $V$.  The $G$-invariance of the
character of $V$ implies that $\res_H V\cong\res_H(aV)$.  Viewing
$V$ and $aV$ as $H$-invariant subspaces of $W$, one can connect
$V$ and $aV$ through a continuous family of $H$-invariant
subspaces along the arc of $S(\rho)$ joining $1$ and $e^{2\pi
i/n}$, in other words, for each $z$ in the arc one can find an
$H$-invariant subspace in the fiber of $\underline{W}$ at $z$ so
that the family of $H$-invariant subspaces varies continuously on
the points $z$.  This is always possible because the set of such
$H$-invariant subspaces of $W$ is homeomorphic to a product of
Grassmann manifolds which is arcwise connected. Translating the
family of $H$-invariant subspaces by the action of $a$ repeatedly
yields the desired $G$-subbundle of $\underline{W}$.

\emph{Case 3:} The case where $\rho(G)=D_n$.  By
Proposition~\ref{prop:extension}, there exist $G_1$- and
$G_\mu$-extensions $V_1$ and $V_\mu$ of $V$ respectively.  Set
$W=\ind_{G_1}^G V_1\oplus\ind_{G_\mu}^G V_\mu$ and consider the
product $G$-vector bundle $\underline{W}$.  Then $V_1$ and $V_\mu$
are contained as $G_1$- and $G_\mu$-submodules in the fibers of
$\underline{W}$ at $1$ and $\mu$, respectively.  Since $\res_H
V_1\cong\res_H V_\mu$, it is possible to connect $V_1$ and $V_\mu$
through a continuous family of $H$-invariant subspaces along the
arc of $S(\rho)$ joining $1$ and $\mu$ as we did in Case 2 above.
We translate it using the action of $b$ and then using the action
of $a$ repeatedly to obtain the desired $G$-subbundle of
$\underline{W}$.
\end{proof}

\begin{rema}
(1) The proof above shows that any $G_1$-extension (and
$G_\mu$-extension when $\rho(G)=D_n$) of $V$ can be realized as
the fiber at $1$ (and at $\mu$ when $\rho(G)=D_n$) of a $G$-vector
bundle over $S(\rho)$.

(2) The proof above also works in the real category except the
extension problem of $V$.  Namely, since
Proposition~\ref{prop:extension} does not hold in the real
category, we need to assume that, in addition to the
$G$-invariance of the character of $V$, $V$ has a $G_1$-extension
when $\rho(G)=O(2)$, and both $G_1$- and $G_\mu$-extensions when
$\rho(G)=D_n$.
\end{rema}

\section{The semi-group structure on $\Vect_G(S(\rho))$}

In this section we determine the semi-group structure on
$\Vect_G(S(\rho))$.  We begin with a simple case.

\begin{lemm} \label{lemm:structure}
Suppose the $G$-action on $S(\rho)$ is effective, in other words,
$\rho\colon G\to O(2)$ is injective.  Then the semi-group
$\Vect_G(S(\rho))$ is generated by
\begin{enumerate}
\item one trivial $G$-line bundle if $\rho(G)\subset SO(2)$,
\item two trivial $G$-line bundles $L^\pm$ if $\rho(G)=O(2)$, and
\item four $G$-line bundles $L^{\pm\pm}$
with a relation $L^{++}+L^{--}=L^{+-}+L^{-+}$ otherwise.
\end{enumerate}
\end{lemm}

\begin{proof}
Since $\rho$ is injective, $P=\rho^{-1}(SO(2))$ acts freely on
$S(\rho)$; so taking orbit spaces by $P$ gives an isomorphism
\[
\Vect_G(S(\rho))\cong \Vect_{G/P}(S(\rho)/P).
\]
In fact, the inverse is given by pulling back elements in
$\Vect_{G/P}(S(\rho)/P)$ by the quotient map from $S(\rho)$ to
$S(\rho)/P$.  Because of this isomorphism, it suffices to study
the semi-group structure on $\Vect_{G/P}(S(\rho)/P)$.  Note that
$S(\rho)/P$ is again a circle or a point, and that the pullback of
a trivial bundle is again trivial.

(1) The case where $\rho(G)\subset SO(2)$.  In this case $P=G$,
i.e., $G/P$ is the trivial group, so the semi-group
$\Vect_{G/P}(S(\rho)/P)$ is generated by one element, that is the
trivial line bundle, as is well known.  This implies (1) in the
lemma.

(2) The case where $\rho(G)=O(2)$.  In this case $G/P$ is of order
two and $S(\rho)/P$ is a point.  Therefore,
$\Vect_{G/P}(S(\rho)/P)$ is generated by two elements of dimension
one.  This implies (2) in the lemma.

(3) The case where $\rho(G)=D_n$ for some $n$.  In this case
$S(\rho)/P$ is again a circle, $G/P$ is of order two, and the
action of $G/P$ on $S(\rho)/P$ is a reflection.  In the sequel it
suffices to treat the case where $n=1$.  But this case is already
studied in~\cite{Kim94}.  (Kim treats real bundles but the same
argument works for complex bundles.)  The result in~\cite{Kim94}
says that $D_1$-vector bundles over $S(\rho)$ are distinguished by
the fiber $D_1$-modules over the fixed points $\pm 1 \in S(\rho)$,
and that any pair of $D_1$-modules of the same dimension is
realized as the fiber $D_1$-modules at $\pm 1$ of a $D_1$-vector
bundle over $S(\rho)$.  Since $D_1$ is of order two, there are two
one-dimensional $D_1$-modules (one is the trivial one $\C_+$ and
the other is the nontrivial one $\C_-$) and that any $D_1$-module
is a direct sum of them.  Therefore, there are four inequivalent
$D_1$-line bundles $L^{\pm\pm}$, where $L^{\epsilon\delta}$
($\epsilon$ and $\delta$ stand for $+$ or $-$) denotes the
$D_1$-line bundle with the fiber $D_1$-modules $\C_\epsilon$ at
$1$ and $\C_\delta$ at $-1$, and the structure of
$\Vect_{D_1}(S(\rho))$ is as stated in (3).

For the reader's convenience, we shall give the argument
in~\cite{Kim94} briefly.  First we observe that any pair of
one-dimensional $D_1$-modules can be realized as the fiber
$D_1$-modules at $\pm 1$ of a $D_1$-line bundle over $S(\rho)$. In
fact, the trivial (real) line bundle and the (real) Hopf line
bundle have respectively two different $D_1$-liftings to the total
space, and each $D_1$-lifting of the trivial (real) line bundle
has the same fiber $D_1$-modules at $\pm 1$ while that of the
(real) Hopf line bundle has different fiber $D_1$-modules at $\pm
1$.  We consider complexification of them.  Then any pair of
$D_1$-modules of dimension $m$ can be realized as the fiber
$D_1$-modules at $\pm 1$ by taking the Whitney sum of suitable $m$
number of those complexified $D_1$-line bundles.  On the other
hand, the same technique used in the proof of
Theorem~\ref{main:condition_of_fiber_$H$-module} shows that any
$D_1$-vector bundle over $S(\rho)$ decomposes into the Whitney sum
of the above $D_1$-line bundles.
\end{proof}

Theorem~\ref{main:condition_of_fiber_$H$-module} applied with
$G=G_\chi$ and the irreducible $H$-module with character $\chi$
says that the assumption in Lemma~\ref{lemm:hom_argument} is
satisfied when $X=S(\rho)$, so $\Vect_{G_\chi}(S(\rho),\chi)$ has
the same semi-group structure as $\Vect_{G_\chi/H}(S(\rho))$. Here
the action of $G_\chi/H$ on $S(\rho)$ is effective, so the lemma
above can be applied to $\Vect_{G_\chi/H}(S(\rho))$.  In the
sequel the semi-group structure on $\Vect_{G_\chi}(S(\rho),\chi)$
is divided into three types depending on $\rho(G_\chi)$.  We
denote the generators of $\Vect_{G_\chi}(S(\rho),\chi)$
corresponding to the generators in Lemma~\ref{lemm:structure} by
\[
\begin{cases}
L_\chi, \quad&\text{if $\rho(G_\chi)\subset SO(2)$,}\\
L_\chi^{\pm}, \quad&\text{if $\rho(G_\chi)=O(2)$,}\\
L_\chi^{\pm\pm}, \quad&\text{otherwise.}
\end{cases}
\]

The following theorem follows immediately from
Lemma~\ref{lemm:structure}.

\begin{main}
The semi-group $\Vect_{G_\chi}(S(\rho),\chi)$ is generated by
\begin{enumerate}
\item one element $L_\chi$ if $\rho(G_\chi)\subset SO(2)$,
\item two elements $L_\chi^\pm$ if $\rho(G_\chi)=O(2)$, and
\item four elements $L_\chi^{\pm\pm}$ with the relation $L_\chi^{++}+
L_\chi^{--}=L_\chi^{+-}+L_\chi^{-+}$ otherwise.
\qed
\end{enumerate}
\end{main}

\begin{rema}
If $\rho(G)\subset SO(2)$, then $\rho(G_\chi)$ is of type (1)
above for any $\chi$.  If $\rho(G)=O(2)$, then $\rho(G_\chi)$ is
of type (1) or (2) above; more precisely $\rho(G_{\chi})=SO(2)$ or
$O(2)$ because the $G$-action on $\Irr(H)$ reduces to an action of
$G/H=\rho(G)=O(2)$ and the action of $SO(2)$ on $\Irr(H)$ is
trivial since $SO(2)$ is connected.  Moreover, if $\rho(G)=D_n$,
then $\rho(G_\chi)$ is of type (1) or (3) above.  Therefore, the
semi-group structure on $\Vect_G(S(\rho))$ can be read from the
theorem above and the isomorphism $\Phi$ in Section 2.
\end{rema}

Using the theorem above, one can easily enumerate $G$-vector
bundles over $S(\rho)$ with a fixed $H$-module $V$ as the fiber
$H$-module.  Since $V$ must have a $G$-invariant character by
Theorem~\ref{main:condition_of_fiber_$H$-module}, one can express
the character of $V$ as
\[
\sum_{\chi\in\Irr(H)/G} m_\chi \Bigl(\sum_{\lambda\in
G(\chi)}\lambda\Bigr)
\]
with non-negative integers $m_\chi$, where $m_\chi$'s are zero for
all but finitely many $\chi$'s in $\Irr(H)/G$ because $V$ is of
finite dimension.  Set
\[
e(\chi)=
\begin{cases} 0, \quad&\text{if $\rho(G_\chi)\subset SO(2)$}\\
1, \quad&\text{if $\rho(G_\chi)=O(2)$}\\
2, \quad&\text{otherwise.}
\end{cases}
\]

With this understood

\begin{coro} \label{coro:enumeration}
The number of isomorphism classes of $G$-vector bundles over
$S(\rho)$ with $V$ as the fiber $H$-modules is given by
$\prod_{\chi\in\Irr(H)/G}(m_\chi+1)^{e(\chi)}$.  \qed
\end{coro}

\section{Isomorphism theorem}

In this section we present another approach to study the
semi-group structure on $\Vect_G(S(\rho))$.  The following
theorem, which we call an isomorphism theorem, reduces the study
of $\Vect_G(S(\rho))$ to representation theory.

\begin{theo} \label{theo:isomorphism}
Two $G$-vector bundles $E$ and $E'$ over $S(\rho)$ are isomorphic
if and only if the fiber $G_z$-modules $E_z$ and $E'_z$ at $z\in
S(\rho)$ are isomorphic for $z=1$ (and for $z=\mu$ when
$\rho(G)=D_n$).
\end{theo}

\begin{proof}
The necessity part is obvious, so we prove the sufficiency. We
note that if there exists an equivariant isomorphism $\Psi\colon
E\to E'$, then it must satisfy the equivariance condition
\[
\Psi_{\rho(g)z} = g \Psi_z g^{-1}
\]
for any $g\in G$ where $\Psi_z = \Psi|_{E_z}$.  By the assumption
we have a $G_1$-linear isomorphism $\Psi_1$ (and a
$G_{\mu}$-linear isomorphism $\Psi_\mu$ when $\rho(G)=D_n$).  In
the following we will define $\Psi_z$ for all $z\in S(\rho)$ using
the above equivariance condition to get an equivariant isomorphism
$\Psi$.  We consider three cases according to the images of $G$ by
$\rho$.

\emph{Case~1:} The case where $\rho(G)=SO(2)$ or $O(2)$.  In this
case the $G$-action on $S(\rho)$ is transitive, so for any $z\in
S(\rho)$ we define $\Psi_z=g\Psi_1g^{-1}$ with $g\in G$ such that
$z=\rho(g)1$.  The well-definedness follows from the
$G_1$-equivariance of $\Psi_1$.  This gives the desired
equivariant isomorphism $\Psi$.

\emph{Case~2:} The case where $\rho(G)=Z_n$.  Let $\zeta=e^{2\pi
i/n}$ and define $\Psi_{\zeta}$ by $a\Psi_1a^{-1}$.  The map
$\Psi_\zeta$ is also an $H$-equivariant isomorphism.  We connect
$\Psi_1$ and $\Psi_\zeta$ along the arc of $S(\rho)$ joining $1$
and $\zeta$, in other words, we find an $H$-equivariant linear
isomorphism $\Psi_z$ for each $z$ in the arc of $S(\rho)$ so that
$\Psi_z$ is continuous at those $z$. (This is always possible
because the set of $H$-linear isomorphisms between $E_z$ and
$E'_z$ is arcwise connected, in fact, homeomorphic to a product of
$\GL(N,\C)$'s.) Now we define $\Psi_z$ for any $z\in S(\rho)$
using the equivariance condition $\Psi_{\zeta z}=a\Psi_za^{-1}$.
This gives the desired isomorphism $\Psi$.

\emph{Case~3:} The case where $\rho(G)=D_n$.  Note that the
equivariance condition of $\Psi$ is
\[
\Psi_z = h \Psi_z h^{-1} \text{ for any $h\in H$}, \quad \Psi_{\zeta z} = a
\Psi_z a^{-1}, \quad \Psi_{\bar z} = b \Psi_z b^{-1}.
\]
We connect $\Psi_1$ and $\Psi_\mu$ along the arc joining $1$ and
$\mu$ to obtain $\Psi_z$ for $z$ in the arc.  Then using the
equivariance condition $\Psi_{\bar z}=b\Psi_zb^{-1}$, we define
$\Psi_z$ for $z$ in the arc joining $1$ and $\mu^{-1}$.  Thus we
have defined $\Psi_z$ for $z$ in the arc joining $\mu^{-1}$ and
$\mu$.  We then define $\Psi_z$ for all $z$ using the equivariance
condition $\Psi_{\zeta z}=a\Psi_za^{-1}$.  This gives the desired
isomorphism.
\end{proof}

For a group $K$ we denote by $\Rep(K)$ the set of isomorphism classes of
$K$-modules, and by $\Rep(G_1,G_\mu)$ the set of elements
$(V,W)\in \Rep(G_1)\times\Rep(G_\mu)$ with $\res_HV=\res_HW$.
Restriction of a $G$-vector bundle over $S(\rho)$ to fibers
at $1$ (and $\mu$ when $\rho(G)=D_n$) yields a map
\[
\Gamma\colon \Vect_G(S(\rho))\to
\begin{cases} \Rep(H)^G \quad&\text{if $\rho(G)\subset SO(2)$,}\\
\Rep(G_1)\quad&\text{if $\rho(G)=O(2)$,}\\
\Rep(G_1,G_\mu)\quad&\text{if $\rho(G)=D_n$}
\end{cases}
\]
where $\Rep(H)^G$ denotes the subset of $\Rep(H)$ with
$G$-invariant character. The target of the map $\Gamma$ is a
semi-group under direct sum.  With this understood

\begin{prop} \label{prop:gamma}
The map $\Gamma$ is an isomorphism.
\end{prop}

\begin{proof}
It is obvious that $\Gamma$ is a homomorphism and that the
characters of all representations in $\Rep(G_1)$ (when
$\rho(G)=O(2)$) and $\Rep(G_1,G_\mu)$ (when $\rho(G)=D_n$) are
$G$-invariant. The surjectivity follows from
Theorem~\ref{main:condition_of_fiber_$H$-module} (when
$\rho(G)\subset SO(2)$) and the first remark following it in
Section 4 (when $\rho(G)=O(2)$ or $D_n$), and the injectivity
follows from Theorem~\ref{theo:isomorphism}.
\end{proof}

As a matter of fact, the source and target of the map $\Gamma$
have more structures, that is, they have products given by tensor
product and $R(G)$ acts on them naturally through the tensor
product. In fact, $R(G)$ acts on the target through the
restriction to $H$, $G_1$, or $G_\mu$.  Clearly the map $\Gamma$
preserves these structures.

\section{Triviality of $G$-vector bundles over a circle}

In this section we investigate when a $G$-vector bundle over
$S(\rho)$ is trivial.  Here is the criterion of triviality of a
$G$-vector bundle over $S(\rho)$.

\begin{lemm} \label{lemm:criterion_for_triviality}
\textup{(1)} A $G$-vector bundle over $S(\rho)$ is trivial if and
only if the fiber $G_z$-module at $z=1$ (and at $z=\mu$ when
$\rho(G)=D_n$) extends to a (same when $\rho(G)=D_n$) $G$-module.

\textup{(2)} Unless $\rho(G)\subset SO(2)$, the number of the
isomorphism classes of trivial $G$-vector bundles over $S(\rho)$
with an irreducible fiber $H$-module $V$ agrees with the number of
$G$-extensions of $V$.
\end{lemm}

\begin{proof}
(1) The necessity is trivial and the sufficiency follows from
Theorem~\ref{theo:isomorphism}.

(2) Let $W$ and $W'$ be two $G$-extensions of $V$ and suppose that
the product bundles $\underline{W}$ and $\underline{W'}$ are
isomorphic.  Then $\res_{G_1}W\cong \res_{G_1}W'$ (and
$\res_{G_\mu}W\cong\res_{G_\mu}W'$ if $\rho(G)=D_n$).  It is easy
to see from Proposition~\ref{prop:extension} that each
$G$-extension of $V$ is distinguished by its restriction to $G_1$
(and $G_\mu$ if $\rho(G)=D_n$).  Therefore, $W$ and $W'$ are
isomorphic as $G$-modules, proving (2).
\end{proof}

\begin{main}
The triviality of the generators appeared in
Theorem~\ref{main:semi_group_structure} is as follows.
\begin{enumerate}
\item $L_\chi$ is trivial.
\item $L_\chi^{\pm}$ are both trivial or both nontrivial.
\item Two of $L_\chi^{\pm\pm}$ are trivial and the other two are
nontrivial if $|\rho(G_\chi)|/2$ is odd, and $L_\chi^{\pm\pm}$ are
all trivial or all nontrivial if $|\rho(G_\chi)|/2$ is even.
\end{enumerate}
\end{main}

\begin{proof}
This follows from Proposition~\ref{prop:extension} and
Lemma~\ref{lemm:criterion_for_triviality}.
\end{proof}

In fact, $L_\chi^{\pm}$ are related through the tensor product
with a one-dimensional nontrivial representation $G_\chi\to
O(2)/SO(2)=\{\pm 1\}$ (which is $\rho$ composed with the
projection $O(2)\to O(2)/SO(2)$), and $L_\chi^{\pm\pm}$ are
related through the tensor product with four $G_\chi$-line bundles
$L^{\pm\pm}$ which are pullback of the four $D_1$-line bundles by
the quotient map from $S(\rho)$ to $S(\rho)/P$ where
$P=\rho^{-1}(SO(2))$ as before. Note that the fiber $H$-modules of
$L_\chi^{\pm\pm}$ are trivial. The $G_\chi$-line bundles
$L^{\pm\pm}$ are well understood and one (actually two) of
$L_\chi^{\pm\pm}$ is trivial in case $n$ is odd, so we completely
understand $L_\chi^{\pm\pm}$ in this case.  But we do not know
$L_\chi^{\pm\pm}$ explicitly when $n$ is even and they are all
nontrivial.

\begin{coro} \label{coro:triviality_chi}
If $\rho(G_\chi)\subset SO(2)$, then every element in
$\Vect_{G_\chi}(S(\rho),\chi)$ is trivial.
\end{coro}

\begin{proof}
The corollary follows from Theorem~\ref{main:semi_group_structure}
(1) and Theorem~\ref{main:triviality_for_L} (1).
\end{proof}

\begin{mcoro}
Every $G$-vector bundle over $S(\rho)$ is trivial if
$\rho(G)\subset SO(2)$.
\end{mcoro}

\begin{proof}
If $\rho(G)\subset SO(2)$, then $\rho(G_\chi)\subset SO(2)$ for
any $\chi$.  Therefore the corollary follows from
Lemma~\ref{lemm:induced_decomposition} and
Corollary~\ref{coro:triviality_chi}.
\end{proof}

\section{Description of $G$-line bundles over a circle}

In this section, we describe $L_\chi^{\pm\pm}$ explicitly when
$\chi$ is the character of a one-dimensional $H$-module.  In the
following, $\varphi$ denotes a one-dimensional $H$-representation
with $G$-invariant character $\chi$.  (Actually $\chi$ agrees with
$\varphi$ since $\varphi$ is one-dimensional.)  To simplify
notation we denote $G_\chi$ by $G$.

When $\rho(G)=D_n$, let $a$ and $b$ be as before, i.e., they
denote elements of $G$ whose images by $\rho$ are respectively the
rotation through an angle $2\pi/n$ and the reflection by $x$-axis.
Note that $G$ is generated by $H$, $a$ and $b$ under the relations
$a^n$, $b^2$, and $(ab)^2\in H$.

\begin{lemm} \label{lemm:line}
Suppose $\rho(G)=D_n$.  Then
\begin{enumerate}
\item $\varphi(a^n)^2=\varphi(abab^{-1})^n$,
\item when $n$ is even, $\varphi$ has a $G$-extension if and
only if $\varphi(a^n)=\varphi(abab^{-1})^{n/2}$.
\end{enumerate}
\end{lemm}

\begin{rema}
When $n$ is odd, we know that $\varphi$ has a $G$-extension by
Proposition~\ref{prop:extension} (3).  It can also be seen from
the proof of (2) below.
\end{rema}

\begin{proof}
Let $\widetilde\varphi\colon G_1\to \GL(1,\C)$ be an extension of
$\varphi$ to $G_1$.  Since (the character of) $\varphi$ is
$G$-invariant (in particular, $G_1$-invariant) and $H$ is an index
two subgroup of $G_1$, such an extension exists by
Proposition~\ref{prop:extension} (1).

(1) It is elementary to see that $a^iba^i\in G_1$ and
$a^iba^iaba\in H$ for $1\le i\le n-1$.  Since $\varphi$ is
$G$-invariant we have
\[
\widetilde\varphi(a^iba^i)\widetilde\varphi(aba)=\varphi(a^iba^iaba)
=\varphi(a(a^iba^iaba)a^{-1})=\widetilde\varphi(a^{i+1}ba^{i+1})
\widetilde\varphi(b)
\]
for $1\leq i\leq n-1$.  By an inductive application of the above
identity, we have
\[
\widetilde\varphi(aba)^n=\widetilde\varphi(a^nba^n)
\widetilde\varphi(b)^{n-1}.
\]
It follows that
\[
\varphi(abab^{-1})^n=\widetilde\varphi(aba)^n
\widetilde\varphi(b^{-1})^n
=\widetilde\varphi(a^nba^n)\widetilde\varphi(b)^{n-1}\widetilde
\varphi(b^{-1})^n
=\varphi(a^n)^2.
\]

(2) The necessity is obvious, so we shall prove the sufficiency.
Let $A$ be an $n$-th root of $\varphi(a^n)$.  Then $(A^{-2}
\varphi(abab^{-1}))^n=1$ by the identity in (1) above, so there is
an integer $k$ (determined module $n$) such that
$A^{-2}\varphi(abab^{-1})=\zeta^k$ where $\zeta=e^{2\pi i /n}$.
The equality $\varphi(a^n)=\varphi(abab^{-1})^{n/2}$ is equivalent
to $k$ being even.  We define $\widetilde\varphi(a)=A\zeta^{k/2}$.
Then $\widetilde\varphi(a)^n =(A\zeta^{k/2})^n=A^n=\varphi(a^n)$.
Therefore, to see that the extended $\widetilde\varphi$ is a
$G$-extension of $\varphi$, it only remains to check that
$\widetilde\varphi(a)^2\widetilde\varphi (b)^2=\varphi((ab)^2)$.
(Remember that $(ab)^2\in H$.) The left hand side at the identity
is equal to $A^2\zeta^k\widetilde\varphi(b)^2$ while the right
hand side is equal to $\varphi(abab^{-1})\varphi(b^2)$ which
agrees with $A^2\zeta^k\widetilde\varphi(b)^2$ because
$A^{-2}\varphi(abab^{-1})=\zeta^k$ by the choice of $k$ above.
\end{proof}

Let $\varphi\colon H\to U(1)$ be a unitary representation with
$G$-invariant character.  If $\rho(G)=D_n$, then there are exactly
four $G$-line bundles with $\varphi$ as the fiber
$H$-representation by Theorem~\ref{main:semi_group_structure} (3).
They are described explicitly in the following example.

\begin{exam} \label{exam:line}
Assume that $\rho(G)=D_n$.  Let $\varphi\colon H\to U(1)$ be a
unitary representation with $G$-invariant character, and let
$\widetilde\varphi\colon G_1\to U(1)$ be a $G_1$-extension of
$\varphi$.  Let $A$ be an $n$-th root of $\varphi(a^n)$.  As
observed in the proof of Lemma~\ref{lemm:line} (2), there is an
integer $k$ such that
$A^{-2}\widetilde\varphi(abab^{-1})=\zeta^k$.  One can check that
\[
h(z,v)=(z,\varphi(h)v)\text{ for $h\in H$}, \quad a(z,v)=(\zeta z, Av), \quad
\text{and}
\quad b(z,v)=(\bar z,\widetilde\varphi(b)z^kv)
\]
define an action of $G$ on $S^1\times\C$.  In fact, it defines a
$G$-line bundle over $S(\rho)$ such that the fiber representation
at $1$ is $\widetilde\varphi$ and that at $\mu$ is given by
\[
h \mapsto \varphi(h), \quad ab \mapsto
A\mu^k\widetilde\varphi(b).
\]
Since the integer $k$ is only determined modulo $n$ (once $A$ is
chosen) and $\mu^n=-1$, this construction gives two $G$-line
bundles with $\widetilde\varphi$ as the fiber representation at
$1$.  (If we take $k+n$ instead of $k$, then the fiber
representation at $\mu$ evaluated on $ab$ changes the sign.) Since
there are two $G_1$-extensions of $\varphi$ by
Proposition~\ref{prop:extension}, the above construction describes
all the four $G$-line bundles over $S(\rho)$ with $\varphi$ as the
fiber $H$-module.
\end{exam}

We now have the classification result for $G$-line bundles over
$S(\rho)$.

\begin{theo} \label{theo:line}
Let $\varphi\colon H\to U(1)$ be an $H$-representation with
$G$-invariant character.  Let $N$ be the number of $G$-line
bundles over $S(\rho)$ with $\varphi$ as the fiber $H$-module.
\begin{enumerate}
\item If $\rho(G)\subset SO(2)$, then $N=1$ and the bundle is
trivial.
\item If $\rho(G)=O(2)$, then $N=2$ and both bundles are
trivial or both are nontrivial.
\item If $\rho(G)=D_n$, then $N=4$ and all the four bundles
are given in Example~\ref{exam:line}, and
\begin{enumerate}
\item if $n$ is odd, then two of them are trivial
and the other two are nontrivial;
\item if $n$ is even, then all the four bundles are
trivial if $\varphi(a^n)=\varphi(abab^{-1})^\frac{n}{2}$, and all
are nontrivial otherwise.
\end{enumerate}
\end{enumerate}
\end{theo}

\begin{proof}
This follows from Theorem~\ref{main:semi_group_structure},
Theorem~\ref{main:triviality_for_L}, and the observation done in
Example~\ref{exam:line}.
\end{proof}

\section{The case when $G$ is abelian}

When $G$ is abelian (and hence so is the subgroup $H$), any
irreducible $H$-module is one-dimensional and $G_\chi=G$ for any
character $\chi$ of $H$.  Therefore $\Vect_{G}(S(\rho))$ is
generated by $G$-line bundles.  Moreover, since $\rho(G)$ is an
abelian subgroup of $O(2)$, it is contained in $SO(2)$ or
isomorphic to $D_1$ or $D_2$.  When $\rho(G)\subset SO(2)$, any
$G$-line bundle over $S(\rho)$ is trivial as we know.  When
$\rho(G)=D_2$, the condition
$\varphi(a^n)=\varphi(abab^{-1})^{n/2}$ in
Theorem~\ref{theo:line}~(3--ii) for $n=2$ holds because $G$ is
abelian; so any $G$-line bundle over $S(\rho)$ is trivial in this
case, too.  But there are two nontrivial $G$-line bundles when
$\rho(G)=D_1$ as claimed in Theorem~\ref{theo:line}~(3--i).  The
following example is simply an interpretation of
Example~\ref{exam:line} to the special case when $\rho(G)=D_1$.

\begin{exam} \label{exam:line_for_abelian}
Suppose $G$ is abelian and $\rho(G)=D_1$.  Then $G=G_1$.  Let
$\varphi\colon H\to U(1)$ be a unitary representation of $H$. (Any
$\varphi$ is $G$-invariant because $G$ is abelian.) As in
Example~\ref{exam:line} choosing a $G$-extension
$\widetilde\varphi$ of $\varphi$ induces a $G$-action on
$S(\rho)\times\C$ defined by
\[
h(z,v)=(z,\varphi(h)v) \text{ for $h\in H$} \quad \text{and}\quad
b(z,v)=(\bar z,\widetilde\varphi(b)zv).
\]
It gives a nontrivial $G$-line bundle over $S(\rho)$. Since there
are two $G$-extensions of $\varphi$, this produces two nontrivial
$G$-line bundles over $S(\rho)$ with $\varphi$ as the fiber
$H$-module.
\end{exam}

Summing up, we have

\begin{prop} \label{prop:line_for_abelian}
Suppose $G$ is abelian and let $E$ be a $G$-vector bundle over
$S(\rho)$.
\begin{enumerate}
\item If $\rho(G)\neq D_1$, then $E$ is trivial.
\item If $\rho(G)=D_1$, then $E$ is the Whitney sum of trivial
bundles and the nontrivial line bundles in
Example~\ref{exam:line_for_abelian}.
\end{enumerate}
\end{prop}

\section{Equivariant $K$-groups of a circle}

In this section we apply the results discussed in the previous
sections to the calculation of the reduced equivariant $K$-group
of $S(\rho)$.  For a compact $G$-space $X$, the equivariant
$K$-group $K_G(X)$ of $X$ is defined to be the Grothendieck group
of finite dimensional $G$-vector bundles over $X$.  If $X$ has a
base point $*$ fixed by the $G$-action, then the reduced
equivariant $K$-group $\widetilde K_G(X)$ is defined to be the
kernel of the restriction homomorphism $K_G(X)\to K_G(*)$ induced
from the inclusion map.  In fact, $K_G(X)$ and $\widetilde K_G(X)$
are algebras over $R(G)$ (although there is no identity element in
$\widetilde K_G(X)$).

The additive structure on $K_G(S(\rho))$ can be determined
completely by Lemma~\ref{lemm:induced_decomposition} and
Theorem~\ref{main:semi_group_structure}. One can also describe the
$R(G)$-algebra structure in terms of representation ring through
the map $\Gamma$ in Section 6.  In the following, we shall compute
$\widetilde K_G(S(\rho))$.  Note that $\widetilde K_G(S(\rho))$ is
defined only when $S(\rho)$ has a fixed point, i.e., $G=H$ or
$\rho(G)=D_1$, and that $\widetilde K_G(S(\rho))$ is trivial if
$G=H$.  Suppose $\rho(G)=D_1$.  Then the $G$-fixed point set
$S(\rho)^G$ consists of two points $\{\pm1\}$ and we take $-1$ to
be a base point.  It follows from Theorem~\ref{theo:isomorphism}
that the restriction homomorphism
\[
\widetilde K_G(S(\rho)) \to \widetilde K_G(S(\rho)^G)\cong R(G)
\]
to fibers at $1$ is injective.
The following theorem determines the image of the
homomorphism as an ideal of $R(G)$, which extends Y.~Yang's result
for $G$ finite cyclic~\cite[Theorem A]{Yan95} to any compact Lie group $G$.
Denote by $\C_+$ and $\C_-$
the $G$-modules of dimension one induced from the trivial and the
nontrivial $D_1$-modules of dimension one, respectively, by the
homomorphism $G\to G/H\cong D_1$.

\begin{theo}
If $\rho(G)=D_1$, then $\widetilde K_G(S(\rho))$ is isomorphic to
the ideal $R(G)(\C_+-\C_-)$ in $R(G)$ generated by $\C_+-\C_-$. In
particular, $\widetilde K_G(S(\rho))$ is torsion-free for any
compact Lie group $G$.
\end{theo}

\begin{proof}
The remark (1) at the end of Section~4 implies that
$R(G)(\C_+-\C_-)$ is contained in the image of $\widetilde
K_G(S(\rho)) \to \widetilde K_G(S(\rho)^G)$, so we prove the converse.

Choose an element $E-F$ in $\widetilde K_G(S(\rho))$.  Then the
fibers of $E$ and $F$ at the base point $-1$ are isomorphic as
$G$-modules.  In particular, $E$ and $F$ have the same fiber
$H$-module and thus $\res_H E_1$ is isomorphic to $\res_H F_1$.
Hence, one can express the image of $E-F$, i.e., $E_1-F_1$ in
$R(G)$ as $E_1-F_1\cong\bigoplus (E_i-F_i)$ where $E_i$ and $F_i$
are irreducible $G$-submodules of $E_1$ and $F_1$, respectively,
such that $\res_H E_i\cong\res_H F_i$.  Here we note that an
irreducible $G$-module $W$ is uniquely determined by $\res_H W$ if
it is reducible, because $2W\cong \ind_H^G\res_H W$ in this case,
see {\cite[Theorem 7.3 (ii), Chapter VI]{BrtD85}}.  Therefore, if
$\res_H E_i\cong \res_H F_i$ is reducible, then $E_i-F_i=0$ in
$R(G)$.  On the other hand, if $\res_H E_i\cong\res_H F_i$ is
irreducible, then both $E_i$ and $F_i$ are $G$-extensions of
$\res_H E_i$.  Thus $F_i$ is isomorphic to $E_i$ or
$E_i\otimes\C_-$ by the last statement of
Proposition~\ref{prop:extension}.  It follows that $E_i-F_i$ is
either zero or $E_i\otimes(\C_+-\C_-)$ in $R(G)$.  Therefore, the
image of $E-F$ is contained in the ideal $R(G)(\C_+-\C_-)$.
\end{proof}

\providecommand{\bysame}{\leavevmode\hbox to3em{\hrulefill}\thinspace}

\end{document}